\documentclass[reqno]{amsart}
\usepackage{epsfig,graphicx}

\def\E{\ifmmode{\mathbb E}\else{$\mathbb E$}\fi} 
\def\N{\ifmmode{\mathbb N}\else{$\mathbb N$}\fi} 
\def\R{\ifmmode{\mathbb R}\else{$\mathbb R$}\fi} 
\def\Q{\ifmmode{\mathbb Q}\else{$\mathbb Q$}\fi} 
\def\C{\ifmmode{\mathbb C}\else{$\mathbb C$}\fi} 
\def\H{\ifmmode{\mathbb H}\else{$\mathbb H$}\fi} 
\def\Z{\ifmmode{\mathbb Z}\else{$\mathbb Z$}\fi} 
\def\P{\ifmmode{\mathbb P}\else{$\mathbb P$}\fi} 
\def\T{\ifmmode{\mathbb T}\else{$\mathbb T$}\fi} 
\def\SS{\ifmmode{\mathbb S}\else{$\mathbb S$}\fi} 
\def\DD{\ifmmode{\mathbb D}\else{$\mathbb D$}\fi} 

\newcommand{\e}{\varepsilon}

\newcommand{\ben}{\begin{enumerate}}
\newcommand{\een}{\end{enumerate}}
\newcommand{\be}{\begin{equation}}
\newcommand{\ee}{\end{equation}}
\newcommand{\bea}{\begin{eqnarray}}
\newcommand{\eea}{\end{eqnarray}}
\newcommand{\bc}{\begin{center}}
\newcommand{\ec}{\end{center}}

\newtheorem{thm}{Theorem}[section]
\newtheorem{cor}[thm]{Corollary}
\newtheorem{lem}[thm]{Lemma}
\newtheorem{prop}[thm]{Proposition}

\theoremstyle{definition}
\newtheorem{defn}{Definition}[section]

\theoremstyle{remark}
\newtheorem{rem}{\rm\bfseries{Remark}}[section]

\newtheorem{exm}[rem]{\rm\bfseries{Example}}

\def\R{{\mathbb R}}

\def\E{{\mathbb E}}
\def\Z{{\mathbb Z}}
\def\C{{\mathbb C}}
\def\R{{\mathbb R}}

\def\N{{\mathbb N}}

\def\BB{{\mathcal B}}
\def\SS{{\mathcal S}}
\def\LL{{\mathcal L}}

\def\DD{{\mathcal D}}

\def\MM{{\mathcal M}}

\def\11{{\mathbb I}}

\def\xbar{{\widetilde x}}
\def\ybar{{\widetilde y}}

\def\Jbar{{\widetilde J}}

\def\n{{\noindent}}

\begin{document}

\title{Floer homology and its continuity for non-compact Lagrangian submanifolds}

\author{Yong-Geun Oh}

\thanks{Partially supported by NSF grant \#DMS 9971446}

\address{Korea Institue for Advanced Study, 
207-43 Cheongryangri-dong Dongdaemun-gu,
Seoul 130-012, KOREA \& 
Department of Mathematics,
University of Wisconsin,
Madison, WI ~53706, USA}
\email{oh@math.wisc.edu}

\begin{abstract}
We give a construction of the Floer homology of the pair of
{\it non-compact} Lagrangian submanifolds, which satisfies 
natural continuity property under the Hamiltonian isotopy which
moves the infinity but leaves the intersection set of the pair
compact. This construction uses the concept of Lagrangian 
cobordism and certain singular Lagrangian submanifolds. We apply
this construction to conormal bundles (or varieties) in the
cotangent bundle, and relate it to a conjecture made by 
MacPherson on the intersection theory of the characteristic 
Lagrangian cycles associated to the perverse sheaves
constructible to a complex stratification on the complex
algebraic manifold.
\end{abstract}
\maketitle

\section{Introduction}

Floer [F1] invented the Floer homology $HF(L_0, L_1)$ of the pair $(L_0, L_1)$ of
Lagrangian submanifolds on symplectic manifolds $(P,\omega )$ with suitable
topological restrictions on the pair.  He defined this by considering the
(generalized) Cauchy-Riemann equation
\begin{equation}
\left\{ \begin{array}{l}
{\partial u\over\partial\tau}
+J{\partial u\over\partial t}=0\cr
u(\tau ,0)\in L_0,\;\;u(\tau ,1)\in L_1 
\end{array}
\right.
\end{equation}
for a map $u:\R\times [0,1]\rightarrow P$.

One crucial property of $HF(L_0,L_1)$ for applications to the problems in
symplectic topology, is the invariance property under the Hamiltonian deformations
of the pair.  Floer's original proof [F1] considers the case where
$L_1=\phi^1_H(L_0)$ and $\pi_2 (P,L_0)=\{ e\}$ where $\phi^1_H:P\to P$ is the
time-one map of the Hamiltonian flow of the function $H:P\times [0,1]\to\R$, and
involves some combinatorial study of the changes occurring to the boundary
operators when a (generic) degenerate intersection occurs between the pairs during
the deformations.  Using the fact that generic types of such degenerate
intersections are either birth-death or death-birth type, he algebraically
analyzed the change. However this study involves a gluing theory of trajectories
on degenerate intersections. Although such a gluing theory is believed to be
possible by now, details were only sketched in [F1].

Because Floer's analysis in [F1] also uses the fact that the action functional is
single-valued in his case (where $\pi_2(P,L_0)=\{ e\}$ is assumed), it was not
clear to the author at the time of writing [O1]
whether this approach can be generalized to more
general cases where the action functional is not single-valued. More
importantly, Floer's original proof does not give naturality of
the chain map. Motivated by Floer's  
approach taken in [F3] for Hamiltonian diffeomorphisms, 
the present author [O1]  used a variant of (1),
\begin{equation}
\left\{ \begin{array}{l}
{\partial u\over\partial\tau}+ J {\partial
u\over\partial t}=0\cr 
u(\tau ,0)\in L, \; u(\tau ,1)\in L_{\rho (\tau)}
\end{array}
\right. 
\end{equation}
for the construction of the chain homomorphim from $HF(L,L_0)$ to $HF(L,L_1)$,
where $\rho :\R\to [0,1]$ is a {\it monotonically increasing}
function with $\rho (-\infty )=0$ and $\rho (+\infty )=1$.  Similar constructions
have been also used in our more recent papers [O2,KO1,2] in relation to a
quantization program of the classical homology theory. In these works,
naturality of the chain map is essential for the analysis of change
of actions and for the continuity proof of symplectic invariants
constructed therein (see [O2] for details). Another way of defining the
chain homomorphism is to transform (1) into the dynamical version 
\begin{equation}
\left\{\begin{array}{l}
{\partial u\over\partial\tau} +J\Big({\partial u\over\partial t}-X_H(u)\Big)=0
\cr u(\tau ,0)\in L,\;\; u(\tau ,1)\in L_0.
\end{array}
\right.
\end{equation} 
The chain map in
this set-up can be defined by making  the Hamiltonian $H$ depend on
$\tau$-variable as in [F3].

However for the cases in [O2,KO1,2] where we consider a family of {\it conormal
varieties} which are non-compact, or more precisely, where the corresponding
Hamiltonian isotopy is no longer compactly supported, the crucial $C^0$-estimates
for the equation turned out to be available for arbitrary
choice of $\rho$ neither in (2) nor (3). The proof of the $C^0$-estimates works
for the perturbed Cauchy-Riemann equation with some particular type of
perturbations which are {\it either compactly supported} as in the most literature
on the Floer theory {\it or directed in certain particular directions} as in [KO2].
In this sense, the present
author's paper [O2] contains a gap in that he overlooked the failure of
the $C^0$-estimate which is needed
for the proof of continuity of Floer homology $HF(H,S:M)$ under the isotopy
of submanifolds $S \subset M$.  The proof of this
$C^0$-estimate in [Section 3, O2] works  for fixed or compactly supported
perturbations of conormal bundles $\nu^*S$,  and turns out to work only
with a particular choice of the function $\rho$ which should be determined
depending on the solution  $u$, for more general types of perturbation.

One purpose of the present paper is to rectify this gap (see Remark 4.3 (1))
by considering a {\it
suspension\/} of (2).  The relevant geometric suspension of Lagrangian
submanifolds is a quite natural operation in symplectic geometry which has been
used in the literature of symplectic topology (see e.g., [A1, Po]).  After we used this
suspension to construct the chain map, it became quite apparent to us that the
idea of our
construction of the chain map applies to more general circumstances, i.e., to
certain Lagrangian cobordisms in $(P,\omega )$.
However, constructing the natural chain map
$$
h_{\mathcal L} :HF_*(L,L_0)\rightarrow HF(L,L_1)
$$
and extending invariance property of the Floer homology
to the case when $L_0$ and $L_1$ are {\it noncompact} and the
Hamiltonian isotopy ${\mathcal L} = \{L_t\}_{0\leq t \leq 1}$ is
{\it not compactly supported} is the main purpose of the present paper.
Surprisingly, this construction involves the notion of {\it Lagrangian cobordism}
and singular Lagrangian submanifolds of
the type that were used by Kasturirangan and the present author in [KO1,2].
This kind of conormal varieties were introduced by mathematicians in the
micro-local analysis (see [GM], [KaSc] for example).

For the rest of the paper, we will always assume that $(P,\omega)$ is tame:
$(P,\omega)$ is called {\it tame} if there exists a compatible
complex structure $J$ such that the metric $g_J:= \omega(\cdot, J\cdot)$ has
bounded sectional curvature and injectivity radius bounded below from zero. We call
such almost complex structure $J$ {\it tame}. It is easy to see that the set of
tame almost complex structures is contractible if non-empty. We will need a more
restricted class of symplectic manifolds which are {\it Weinstein at infinity}
whose definition is referred to [EG1] or to \S 2 of this paper. The following is
the main theorem whose precise statement will be referred to later sections.

\medskip

\n{\bf Theorem I.}~{\it Let $(P,\omega)$ be Weistein at infinity. Let
$L$ and $\LL =\{L_t\}_{0\leq t\leq 1}$ be a (proper) Lagrangian submanifold
and an isotopy of proper Lagrangian submanifolds satisfying ``suitable''
condition at infinity. Suppose $L\cap L_t$ remain compact for all $t\in [0,1]$.
Then there exists a canonical isomorphism
$$
h_{\mathcal L}: HF(L,L_0) \to HF(L,L_1).
$$}

An immediate consequence of the present construction is the following
intersection theorem of the conormal bundles. A similar intersection
result was previously obtained by Eliashberg and Gromov in the name of
``deformed conormal bundles'' using finite dimensional approach of
generating functions [Theorem 0.3.4.1, EG2].

\medskip

\n {\bf Theorem II.} {\it Let $S_1, \, S_2$ be compact submanifolds of $M$
such that $S_1$ is transverse to $S_2$.
Suppose $\phi$ is a Hamiltonian diffeomrorphism
on $T^*M$ of the types or a composition of them
\smallskip

$(1)$ $\phi$ is obtained by a compactly supported Hamiltonian isotopy, or

$(2)$ it is homogeneous symplectomorphic (at infinity) i.e., it is generated
by the Hamiltonian of the form $(q,p) \mapsto \langle p, X_t(q) \rangle$
such that $S_1$ is transverse to $f_t(S_2)$
for all $t$ where $f_t: M \to M$ is the flow of $X_t$, or

$(3)$ it is a fiberwise translation by $df$ where $f$ is a smooth
function defined on the base $M$.
\smallskip

Then
$$
\# (\nu^*S_1 \cap \phi(\nu^*S_2)) \geq \hbox{\rm rank } H_*(S_1\cap S_2)
$$
provided $\nu^*S_1$ is transverse to $\phi(\nu^*S_2)$.
Here $H_*(S_1\cap S_2)$ is in $\Z$-coefficients in the
oriented case and in $\Z_2$-coefficients in general.}

\medskip

We refer to Theorem 7.2 for a more precise statement concerning the
Floer homology of the pair $(\nu^*S_1,\nu^*S_2)$. 

A special case $S_1=M$ and $S_2=S \subset M$ studied in [Oh2] is of
particular interest in relation to the gap in [Oh2] 
mentioned in the beginning. For this case, the
transversality hypothesis in Theorem II is automatically satisfied. 
This leads to complete construction of the chain map and proof of
its continuity property which in turn fills the gap in the proof
of [Theorem 5.4, Oh2]

\medskip
\n{\bf Corollary [Theorem 5.4, Oh2]}~ {\it Denote by $HF_*(S,J:M)$
the Floer homology between $\nu^*S$ and $o_M (= \nu^*M)$. Let
$S^\alpha$ and $S^\beta$ be two isotopic submanifolds of $M$.
Then there is a canonical isomorphism
$$
h_{\alpha\beta}:HF_*(S^\alpha,J^\alpha:M) 
\to HF_*(S^\beta, J^\beta:M)
$$
that preserves the grading.}

\medskip

Next, we like to compare
the intersection result in Theorem II or Theorem 7.2 with the conjecture stated in [GM],
whose precise meaning ought to be clarified. The results from [KO1,2] and
the present paper can be considered as some steps towards this direction.

\medskip

\n{\bf Conjecture} [GM]. {\it Let $\SS_1$ and $\SS_2$ be two
complex stratifications of a complex manifold $X$. Assume they are transverse
to each other. Let $F_1$ and $F_2$ be perverse sheaves constructible with
respect to $\SS_1$ and $\SS_2$. Let $F_1\otimes F_2$ be the tensor product of
$F_1$ and $F_2$ on $X$. Then the global homology groups $H_i(X; F_1\otimes F_2)$
can be computed as Floer homology of $(-1)_*Ch(\chi F_1)$ and $Ch(\chi F_2)$.}
\medskip

The case considered in Theorem II is a special case of the
Fary functors $F_i$ constructible with respect to  the stratifications
$$
\SS_i = \{ S_i, M - S_i\}
$$
for $i = 1, \, 2$ such that their corresponding constructible functions are given by
$$
\chi F_i = \left\{ \begin{array}{ll}
1 & \, x \in S_i \\
0 & \, x \in M - S_i. 
\end{array}\right. 
$$
One can easily check that the characteristic Lagrangian cycle of
$F_i$ is nothing but $\nu^*S_i$.

Beside the conormal varieties considered in [O2,KO1,2], good examples
to which we can apply the construction of the present paper will be
the symplectic manifolds with contact type boundary and
proper Lagrangian submanifolds in them.  We refer to \S 5 [KhSe]
for some relevant discussions of the latter examples which occur naturally
in the study of vanishing cycles of the singularity of
holomorphic functions. See also Remark 4.3 of the present paper
where our construction is applied to answer some question raised in [Khse]
which concerns naturality of certain isomorphism between the Floer homology of
non-compact Lagrangian submanifolds.
While this paper was in the stage of completion, we learned from
K. Hori (see [HIV]) that some interesting class of non-compact Lagrangian cycles
(``wave front trajectories'' they call), which are closely related
to the vanishing cycles of holomorphic Morse function (``super-potential''),
  play an important role in the
mirror symmetry of open strings in the context of Landau-Ginzburg
model through the Picard-Lefschez theory.

\section{Hamiltonian deformations and $C^0$-estimates}

In this section, we review the usual construction [F2,O1,2] of the chain
map under {\it compactly supported} Hamiltonian isotopies.  Let
$j=\{J_t\}_{0\leq t \leq 1}$ be a family
of almost complex structures that is $t$-independent at infinity, say,
$J_t(x) = J_\infty(x)$ at infinity for some almost complex structure $J_\infty$.
We denote by $\hbox{\rm Supp }j$ to be the subset
$$
\hbox{\rm Supp }j = \cup_{t \in [0,1]} \overline{\{x \in P ~|~ J_t(x)\neq J_\infty(x)\}}.
$$
Let $\{ L_s\}_{0\le s\le 1}$ be a Hamiltonian isotopy
associated to a {\it compactly supported\/} Hamiltonian
$$
H:P\times [0,1]\rightarrow \R.
$$
We choose a cut-off function $\rho :\R \to
[0,1]$ such that

\begin{eqnarray*}
\rho  & = & \left\{  \begin{array}{l}
 0 \quad \mbox{\rm for } \, \tau \le 0
\cr 
 1\quad  \mbox{\rm for } \, \tau\ge 1
 \end{array} \right. \\ 
\rho' & \geq & 0,
\end{eqnarray*}

and $\rho_K(\tau )=\rho \big({\tau\over K}\big)$.  The
construction of Floer's chain map
$$
h:HF(L,L_0)\rightarrow HF(L,L_1)
$$
is given by considering either
\begin{equation}
\left\{ \begin{array}{l}
{\partial u\over\partial\tau}+ J(t,u) {\partial u\over\partial t}=0 \cr
u(\tau ,0)\in L, \; u(\tau ,1)\in L_{\rho_K (\tau )}
\end{array} \right.
\end{equation}  
or
\begin{equation}
\left\{ \begin{array}{l}
{\partial u\over\partial\tau}+J(t,u)\Big({\partial u\over\partial
t}-X_{H^{\rho_K(\tau )}}(u)\Big)=0\cr 
u(\tau ,0)\in L,\;\; u(\tau ,1)\in L_0.\end{array}
\right.
\end{equation}
This construction works as long as $(P,\omega )$ is tame
and the deformation $\{ L_s\}_{0\le s\le 1}$ can be realized by an ambient
Hamiltonian isotopy associated to {\it compactly supported\/} Hamiltonian.

\medskip

Recall from [EG1] that a symplectic manifold $(P,\omega)$ is called
convex at infinity if it carries a vector field $X$ which is {\it complete
symplectically dilating} at infinity: A vector field $X$ is complete
symplectically dilating if the flow $\{X^t\}$ of X is complete and
satisfies $(X^t)^*\omega = e^t \omega$. We assume that $(P,\omega)$ allows
an exhausting pluri-subharmonic funtion at infinity. Following [EG1], we call such
manifold {\it Weinstein} (at infinity). We choose
$\varphi$  an exhausting pluri-subharmonic
function with respect to a tame almost complex structure $J$.
We also assume that $J$ is invariant under the flow of $X$ outside
a compact set. Then the level set $\varphi^{-1}(R)$ for sufficiently large
$R$ carries the induced contact structure (in fact a $CR$-structure) on it.
The following $C^0$-estimate can be proven by a version of strong maximum
principle (See [EHS]).

\medskip

\begin{thm} 
Let $j=\{J_t\}_{0\leq t \leq 1}$ be a family of almost complex structures
such that $J_t = J$ outside a compact set.
Let $H: P\times [0,1] \to \R$ be a compactly supported
Hamiltonian. Suppose that $L_0 \cap L_1$ are compact and $L_i$'s are
transverse to the level sets of $\varphi$ at infinity. Then there exists a compact
subset $K = K(P,\omega, \mbox{\rm supp }j, \varphi) \subset P$ such that
$$
\hbox{\rm Image } u \subset K
$$
for all solutions $u$ of (4) or (5).
\end{thm}

\medskip

\begin{proof} 
Consider the function $\varphi\circ u: \R \times [0,1] \to \R$.
Since this function is subharmonic at infinity with respect to the metric induced by $J$,
it has no interior maximum point outside of $\hbox{\rm Supp } j$.
Suppose that it has a maximum at a boundary point outside $\hbox{\rm Supp }j$,
say, at $(\tau_0,1)$ and that
$$
R_0:=\varphi(u(\tau_0,1)) > \sup_{L_0\cap L_1}\varphi.
$$
By the strong maximum principle, we have
\begin{equation}
{\partial \over \partial t}(\varphi\circ u) (\tau_0,1) > 0
\end{equation}
unless $\varphi\circ u$ is constant, which is not possible if
$R_0 > \sup\{\varphi\circ u(\infty), \varphi\circ u(-\infty)\}$.
We note that
\begin{eqnarray*}
{\partial u \over \partial t}(\tau_0,1) & = & J {\partial u\over \partial \tau}
(\tau_0,1) \in J\cdot TL_1 \\
{\partial  \over \partial t} (\varphi\circ u) & = & d\varphi (J
{\partial u\over \partial \tau })
\end{eqnarray*}
On the other hand, we must have
$$
d\varphi \Big({\partial u\over \partial \tau}\Big)(\tau_0,1) = 0
$$
at the maximum point $(\tau_0,1)$.
This implies
$$
{\partial u\over \partial \tau}(\tau_0,1) \in TL_1 \cap T(\varphi^{-1}(R_0))
$$
where $R_0 = \varphi(u(\tau_0,1))$. Since $TL_1 \cap T(\varphi^{-1}(R_0))$
is Legendrian in $\varphi^{-1}(R_0)$ with respect to the induced
contact structure (in fact, the induced $CR$-structure) by the assumption
that $L_1$ is transverse to $\varphi^{-1}(R)$ for sufficiently large $R$,
$J{\partial u\over \partial \tau}(\tau_0,1)$ is tangent to the contact
distribution, which implies
$$
{\partial \over \partial t}(\varphi\circ u) (\tau_0,1)
= d\varphi\Big({\partial u\over \partial t}\Big)(\tau_0,1) =
d\varphi\Big(J {\partial u\over \partial \tau}\Big)(\tau_0,1) =0
$$
This gives rise to contradiction to (2.3).
\end{proof}

\medskip

Examples of convex symplectic manifolds include cotangent bundles of compact
manifolds. Products of two convex manifolds are also convex. The sum of
(exhausting) pluri-subharmonic functions will provide an (exhausting)
pluri-subharmonic function on the product.

Note that Theorem 2.1 already
takes care of the case when Hamiltonian isotopies are compactly supported.
However in relation to the quantization program illustrated by [O2] and [KO1,2],
one needs to consider certain deformations $\{ L_s\}_{0\le s\le 1}$ of conormal
type which cannot be realized by compactly supported Hamiltonians. For example,
consider an isotopy $\{ S^s\}_{0\le s\le 1}$ of submanifolds $S^s\subset M$.  In
[O2], we consider the corresponding deformation of the conormal bundles
$$
\{ \nu^*S^s\}_{0\le s\le 1}\subset T^*M.
 $$
This deformation is realized by the
Hamiltonian
\begin{equation}
H(q, p,s)=\langle p, X_s(q)\rangle
\end{equation}
where $X_s$ is the vector field
realizing the isotopy $\{ S^s\}$ i.e.\ $X_s={d\over ds}\big| S^s$. Certainly, this
Hamiltonian is not compactly supported.  If one naively attempts to do the similar
construction using (4) or (5) one would immediately encounter a problem in
establishing the $C^0$-estimate. In the next sections, we will carry out
construction of the chain map using ``suspension'' which covers this case as a
special case. In hindsight, to get the required $C^0$ estimates, one has to use
a ``good'' choice of the function $\rho$ in (2) which itself will enter in the
Cauchy-Riemann equation and  should be determined.

\section{Lagrangian cobordism}

In this section, we introduce an equivalence relation on the space of Lagrangian
embeddings in a given symplectic manifold $(P,\omega )$. Compare with [A1,C].
\medskip

\begin{defn} 
We say that two Lagrangian submanifolds $L_0$ and $L_1$
are {\it Lagrangian cobordant\/} on $(P,\omega )$ if there exists a Lagrangian
submanifold
$$
 \beta\subset (P,\omega )\times T^*\R
$$
such that
\smallskip
\item{(i)} $\partial\beta =L_0\times \{ (1,0)\}-L_1\times \{ (0,0)\}$
\smallskip

\item{(ii)} $\beta$ has flat collars near $\partial\beta$, i.e.,
$$
\beta = \left\{
\begin{array}{l}
L_1\times \{ (s,0)\} \quad\mbox{\rm for }\, 0\le s\le\e\cr 
L_0\times \{ (s,0)\} \quad\mbox{\rm for } \, 1-\e\le s\le 1.
\end{array} \right. 
$$
for some $\e >0$. We denote by $L_0
{\displaystyle{\mathrel{\mathop{\sim}_\beta}}} L_1$ if $L_0$ and $L_1$ are {\it
Lagrangian cobordant via\/} $\beta$.
\end{defn}

Note that $P\times T^*\R$ with the obvious product symplectic structure is tame if
$(P,\omega)$ is so.

\medskip
\begin{exm}

\item{(1)} Let $L_1=\phi^1_H(L_0)$ for some Hamiltonian $H:P\times
[0,1]\rightarrow \R$.  We may re-choose $H$ so that $H\equiv 0$ for $t$ near $0$
and $1$.  We define the Lagrangian cobordism
$$
\beta_H\subset P\times  T^*\R
$$
by
$$
\beta_H=\{ (x,s,a)\in P\times  T^*\R\mid x\in L_s,\; a=-H(x,s), 0\le s\le
1\}.
$$
One can easily check that $\beta_H$ is Lagrangian and satisfies both (i)
and (ii). Therefore Hamiltonian isotopies are special cases of Lagrangian cobordism.
\smallskip

\item{(2)} We would like to  separately consider the special case of (1) which
was considered in [O2]. Let $\{ S^s\}_{0\le s\le 1}$ be a smooth family of
submanifolds in a smooth manifold $M$, and $\{ \nu^*S^s\}_{0\le s\le 1}$ be their
conormal bundles. We are given an ambient isotopy $\{ \psi^s\}_{0\le s\le 1}$ such
that
$$
S^s=\psi^s(S^0)
$$
and $\{ X_s\}_{0\le s\le 1}$ is its generating vector
fields, then the corresponding Lagrangian cobordism is given by
$$
\{ (q,p,s,a)\in T^*M\times T^*[0,1]\mid q\in S^s, p\in\nu^*_qS^s, a=-\langle
p,X_s(q)\rangle\;\hbox{and}\; 0\le s\le 1\}.
$$
One can easily check that this
becomes flat if we choose the isotopy to be constant near $s=0$ and $1$.
Furthermore, this bordism itself is nothing but the conormal to the suspension
$$
\{ (q,s)\in M\times [0,1]\mid q\in S^s\}
$$
in $T^*(M\times [0,1])$.
\end{exm}

\section{Construction of chain maps} 

In this section, we attempt to construct the chain map
$$
h_\beta :HF(L,L_0)\rightarrow HF(L,L_1)
$$
when $L_0{\displaystyle{\mathrel{\mathop{\sim}_\beta}}} L_1$.
In the beginning, we do not impose any condition on $L_1$ or $L_0$.
Due to the assumption of
flatness near $\partial\beta$, we can smoothly add to $\beta$ two ends $$
L_0\times (-\infty ,0]\times \{ 0\}\amalg L_1\times [1,\infty )\times \{ 0\}. $$
We again denote the resulting manifold by $\beta$.  This $\beta$ will play a role
as the boundary condition at $t=1$ for the Cauchy-Riemann equation that we will
consider.  We still need the boundary condition at $t=0$ which we now describe.

It turns out that the right choice is the following singular Lagrangian
submanifold
$$
\alpha_L :=L\times Ch(1_{[0,1)}) \subset P\times T^*\R
$$
where $Ch(1_{[0,1)})$
is the characteristic Lagrangian cycle of the characteristic function $1_{[0,1)}$
of $[0,1)$ on $\R$ in the sense of [GK]. In fact, we can prove (see [KO2] for the
general case of standard pairs) that
$$
Ch(1_{[0,1)}) = o_\R \big|_{(0,1)}\amalg \nu^*_-(\partial [0,1])\big |_{0}
\amalg \nu^*_+(\partial [0,1])\big |_{1}
$$
where
\begin{eqnarray*}
\nu^*_-(\partial [0,1])\big |_0 & = &\{ (s,a)\in T^*\R ~|~ s=0, a\ge 0 \} \\
\nu^*_+(\partial [0,1])\big |_1 & = & \{ (s,a)\in T^*\R~ |~ s=1, a\ge 0 \} 
\end{eqnarray*}
is the negative and positive part (with respect to the induced orientation)
of the conormal bundle of $\partial [0,1]$ respectively.

\begin{figure}[htb]
\includegraphics{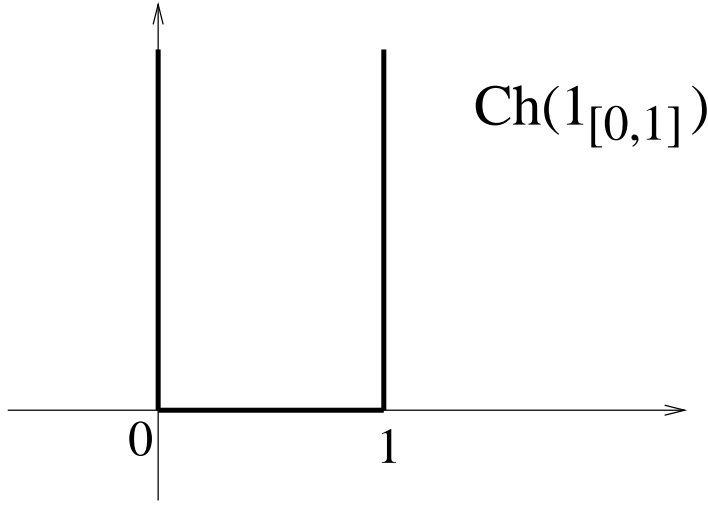}
\caption{$Ch(1_{[0,1)})$}
\end{figure}

In [KO2], we call $Ch(1_{[0,1)})$  the conormal to the {\it standard pair}
$([0,1],\{1\}) $ [GM] and denote it by $\nu^*([0,1],\{1\})$. We refer to [GM]
or [KO2] for the definition of standard pairs.  Then we consider the
following Cauchy-Riemann equation
\begin{equation}
\left\{ 
\begin{array}{l}
{\partial\widetilde u\over\partial\tau}+\widetilde J
\Big({\partial\widetilde u\over\partial t}\Big) =0 
\cr
\widetilde u(\tau,0)\in \alpha_L,\; \widetilde u(\tau ,1)\in\beta
\end{array} \right. 
\end{equation}
where $\Jbar=J\oplus i$, $\widetilde u=(u,v)\subset P\times T^*\R$
where $v = (s,a)$. Since
$\alpha_L$ is singular, we need to desingularize $\alpha_L$ in a suitable way
as in [KO1,2], which we now describe.  We consider
$$ \alpha_\e :=L\times\Upsilon_\e\subset P\times T^*\R $$
where $\Upsilon_\e$ are approximations of $Ch(1_{[0,1)})$ drawn as

\begin{figure}[htb]
\includegraphics{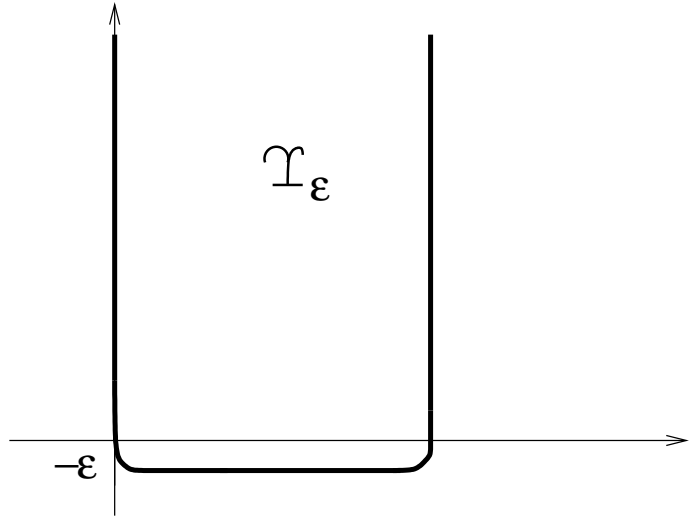}
\caption{Approximation of $Ch(1_{[0,1)})$}
\end{figure}

\n  Since we assume that $\beta$ is flat near $\partial \beta$, we can choose $\e
>0$ so that
$$
\beta\cap\alpha_\e\cap\{ 0<s<\e\;\hbox{ or }\; 1-\e <s<1,\;\hbox{and }\; a=0\}
=\emptyset.
$$
On the other hand, we have
\begin{eqnarray}
\partial\beta\cap\alpha_\e\cap\{ s=0\} & = & L\cap L_0\times \{ (0,0)\}\\
\partial\beta\cap\alpha_\e\cap \{ s=1\}& = & L\cap L_1\times \{(1,0)\}.
\nonumber
\end{eqnarray}
If we assume that $L$ is transverse to both $L_0$ and $L_1$, we can apply
Hamiltonian perturbations in $P\times T^*\R$ of $\beta$ away from the sets (9)
and make $\beta$ intersect transversely with $\alpha_L$. Now for each given $x\in
L\cap L_0$ and $y\in L\cap L_1$, we study the equation
\begin{equation}
\left\{ 
\begin{array}{l}
{\partial\widetilde u\over\partial\tau} +\widetilde
J \Big({\partial\widetilde u\over\partial t}\Big)=0\cr 
\widetilde u(\tau
,0)\in\alpha_L,\; \widetilde u(\tau ,1)\in\beta\cr 
\widetilde u(-\infty
)=\widetilde x=(x,0,0), \, \widetilde u(+\infty)=\widetilde
y=(y,1,0).
\end{array} \right.
\end{equation}

\begin{rem} 
Let us disseminate (10) for the case $\beta = \beta_H$.
In this case, the equation (10) can be re-written as
$$
\left\{ 
\begin{array}{l}
{\partial u\over\partial\tau} +
J {\partial u\over\partial t}=0\cr
u(\tau,0)\in L,\; u(\tau ,1)\in L_{s(\tau,1)} \cr
u(-\infty)=x, \, u(+\infty)= y\cr 
\cr
\overline \partial v = 0 \cr
v(\tau, 0) \in Ch([0,1)), \, a(\tau,1) = -H(u(\tau,1),s(\tau,1)) \cr
s(-\infty) = 0, \, s(\infty)= 1 
\end{array} \right.
$$
The first part of this equation is nothing but (2) with
$\rho(\tau) = s(\tau,1)$ but $s$ itself must be solved.
Furthermore unlike (2), $u$ and $s$ are coupled to each other.
It is rather interesting and mysterious to us that our effort obtaining
the $C^0$-estimates of (2) has led us to considering the coupled
Cauchy-Riemann equation of
$u$ and $\rho$ in the suspended space.
\end{rem}

We note that if $(P,\omega)$ is Weinstein at infinity, so is $(P, \omega) \times
(T^*\R, \omega_0)$. If $\varphi$ is an exhausting pluri-subharmonic function on $P$,
the
$$
\widetilde\varphi(x,s,a) = \varphi(x) + {1\over 2}(s^2 +a^2)
$$
will be an exhausting pluri-subharmonic function on
$(P, \omega) \times(T^*\R, \omega_0)$.
Furthermore both $\alpha_L$ and $\beta$ are {\it fixed}
Lagrangian submanifolds and $\alpha_L$ is transverse to the
level sets of $\widetilde \varphi$ at infinity.
Therefore we have the following a priori $C^0$-estimate
for the solutions of (10) from Corollary 2.2.
\medskip

\begin{prop} 
Suppose that $(P, \omega)$ is Weinstein at infinity.
Assume
\smallskip

 $(1)$ $\alpha_L \cap \beta$ is compact

 $(2)$  $\beta$ is transverse to the level sets of $\widetilde \varphi$ at infinity.
\smallskip

\n Then for any given $\widetilde x,\, \widetilde y \in \alpha_L \cap \beta$,
 there exists a compact subset
$K = K(\widetilde x,\widetilde y, \beta) \subset P \times T^*\R$ such that
$$
\hbox{\rm Image }\widetilde u \subset K
$$
for all $\widetilde u \in \MM_\epsilon (\Jbar, \beta:\widetilde x,\widetilde y)$.
\end{prop}
\smallskip

We now study the moduli-space 
$\MM_\e (\Jbar, \beta)$ 
of solutions of (8) with finite energy. This is decomposed into
$$
\MM_\e (\Jbar, \beta) = \cup_{\widetilde x, \widetilde y \in \alpha_L \cap \beta}
\MM_\e (\Jbar, \beta :\widetilde x,\widetilde y)
$$ 
where $\MM_\e (\Jbar, \beta :\widetilde x,\widetilde y)$ is the set of
solutions of (10). Note that there is a natural $\R$-action  on 
$\MM_\e (\Jbar, \beta :\widetilde x,\widetilde y)$ 
by translations in the $\tau$-direction. We denote 
$$
\widehat{\MM}_\e (\Jbar, \beta :\widetilde x,\widetilde y) 
= \MM_\e (\Jbar, \beta :\widetilde x,\widetilde y)/\R
$$ 
and 
$$
n_\e (x,y :\beta )=\# (\widehat{\MM}_\e (\Jbar, \beta :\xbar,\ybar))
$$
when
$$
\dim \MM_\e (\Jbar, \beta :\xbar,\ybar)=0.
$$
By the standard compactness theorem and the dimension counting arguments, the
zero-dimensional component of $\widehat{\MM_\e}(\Jbar, H, \beta :\xbar,\ybar)$
is compact {\it under suitable assumptions on $L, L_0, L_1$ and $\beta$}.
For example, we
may assume that $L, L_0$ and $L_1$ are monotone in $P$ and $\beta$ is
monotone in $P\times T^*\R$.  (See
[O1]).  We will always assume these conditions from now on for the simplicity
in presenting the main ideas of our construction, although one could
consider more general cases using the sophisticated construction employed in
[FOOO].

We define a map
$$
h_{\beta ,\e}: CF(L,L_0:J,\alpha_\e )\rightarrow CF(L,L_1:J,\alpha_\e )
$$
by
\begin{equation}
h_{\beta ,\e}(x)=\sum_y n_\e (x,y:\beta)y
\end{equation}
and study its chain property.
\medskip

\begin{defn}
Let $(L_0, L_1)$ be a pair of Lagrangian submanifolds
transverse to $L$.  Define $\BB (L_0, L_1)=$ the set of
Lagrangian cobordisms $\beta$ from $L_0$ to $L_1$
$$
\BB_0(L_0, L_1:L)=\{\beta\in\BB (L_0, L_1)\mid \beta \,\,  \mbox{\rm 
is transverse to }\, \alpha_L\}.
$$
\end{defn}

\begin{lem} 
The set $\BB_0(L_0,L_1:L)$ is a residual subset of $\BB
(L_0, L_1)$ in $C^1$-topology.
\end{lem}

\begin{proof}  
It is enough to consider Hamiltonian perturbations of given
$\beta$ that are fixed near $\partial\beta$. The proof of this is standard  which
we omit.
\end{proof}
\medskip

\begin{exm}
Consider a Hamiltonian isotopy from $L_0$ to $L_1$ and its
corresponding Lagrangian cobordism
$$
\beta_H=\{ (x,s,a)\mid x\in L_s,\; a=-H(x,s)\}.
$$
We call this a {\it
Hamiltonian cobordism}. In this case, we note that
$$
\beta_H\cap\alpha_L=\{ (x,s,a)\in P\times T^*\R \mid x\in L\cap
L_s\;\hbox{and}\; a=-H(x,s)=0\}
$$
and $\beta_H$ is transverse to $\alpha_L$ if and only
if
$$
T_xL\oplus T_xL_s=T_xP, \quad {\partial H\over\partial s}(x,s)\not= 0
$$
at each $(x,s,a)\in\beta_H\cap\alpha_L$.  However in general,  we
cannot avoid non-transverse intersections for a one parameter family
$\{ L_s\}_{0\le s\le 1}$, which forces us to look at perturbations of $\beta_H$ on
$P\times T^*\R$ to obtain transversal pairs $(\alpha_L,\beta )$ with $\beta$ close
to $\beta_H$.
\end{exm}

As usual in the Floer theory, we examine
compactness property of the one-dimensional component of
$\widehat{\MM_\e}(\Jbar, \beta :\xbar,\ybar)$ to study
the chain property of $h_{\beta,\e}$, i.e., the identity
\begin{equation}
h_{\beta,\e}\circ \partial_0=\partial_1\circ h_{\beta ,\e}
\end{equation}
\medskip
We consider  one dimensional
components of $\widehat{\MM_\e} (\Jbar, \beta )$ and study structure of the
boundary of each one-dimensional component in its compactification.  Standard
dimension counting argument tells us that the boundary of $\MM_\e (\Jbar, \beta
:\widetilde x,\widetilde z)$ consists of the cusp-trajectories of the form
$\widetilde u_1\#\widetilde
u_2$ where $(\widetilde u_1, \widetilde u_2)$ are elements in $\MM_\e (\widetilde
J, \beta :\widetilde x, \widetilde y)\times \MM_\e(\Jbar, \beta :\widetilde
y,\widetilde z)$.  Here, a priori, $\widetilde y$ could be any element in the
intersection set $\beta\cap\alpha_{L,\e}$, not just in the hypersurface of $s = 0$ or $1$.
 This will prevent us
from associating a chain homomorphism to general Lagrangian cobordism.
From now on, we will mainly concern the case of Hamiltonian cobordism.

\medskip

Let us first examine the condition (1) from Proposition 4.2
 that {\it $\alpha_L \cap \beta$ is
compact}. This is certainly the case if $L$ is compact. For the case of
Hamiltonian cobordism $\beta_H$, it is easy to see that
$\alpha_L\cap \beta_H$ is compact if and only if $L \cap L_t$
is compact for all $t\in [0,1]$. In general, we introduce the following
definition.
\medskip

\begin{defn}
Let $ {\mathcal L} =\{L_t\}_{0\leq t \leq 1}$ be a Hamiltonian
isotopy. We say that {\it intersections do not escape to infinity} if
$\cup_{t \in [0,1]} (L \cap L_t)$ is compact.
\end{defn}

Under this condition, we prove the following proposition,
which will eliminate  those intersections $\tilde y$ away from
$\partial \beta_H \cap \alpha_L$ (i.e., away from $s =0$ or $s=1$)
that provides the obstruction to the existence of chain property.
\medskip

\begin{lem}
Let $L \subset P$ and $ {\mathcal L} = \{L_t\}$ be a
Hamiltonian isotopy of $L_0$ such that the intersections $L\cap L_t$
do not escape to infinity. Let $\beta_H$ be a Hamiltonian cobordism
associated to the Hamiltonian isotopy $\LL$. Then we can change $H$ to $H^\prime$
so that $\phi^t_H = \phi^t_{H^\prime}$, and
\begin{equation}
\beta_{H^\prime} \cap \alpha_L = L \cap L_0 \times \{(0,0)\} \coprod
L \cap L_1 \times \{(1,0)\}
\end{equation}
\end{lem}
\begin{proof} 
We recall
$$
\beta_H \cap \alpha_L =\{ (x,s,a) \in P\times T^*\R~|~
x \in L \cap L_s, \, a = -H(x,s) = 0, \, s\in [0,1]\}
$$
Since $\cup_{s \in [0,1]} L \cap L_s$ is compact by hypothesis,
$\hbox{\rm Image }H|_{\cup_{s \in [0,1]}L \cap L_s}$ is compact.
Therefore we can choose a non-negative function
$$
\chi: [0,1] \to \R_+
$$
so that
\medskip

\item{(i)} $\chi(s) = 0$ for  $s$ near 0 or 1.

\item{(ii)} $\chi(s) + H(x,s) > 0 $ for $(x,s)$ such that
$x \in \cup _{s \in [\delta, 1-\delta]}L \cap L_s$
for some small $\delta >0$.
\medskip

\n We just choose $H^\prime(x,s) := H(x,s) + \chi(s)$
as our new Hamiltonian. 
\end{proof}

\medskip

From now on based on Lemma 4.3 or its proof, we use only the Hamiltonians that   
satisfy
\begin{equation}
H(x,s) > 0 \quad \mbox{\rm for } \,  (x,s) \in \cup_{s \in [\delta,
1-\delta]}L \cap L_s
\end{equation}
for the Hamiltonian cobordism $\beta_H$ when we perform
construction of the chain map $h_{\beta_H,\e}$. 
We call such Hamiltonians {\it (positively) admissible to } $(L, \LL)$.
The following lemma is easy to check
\begin{lem}
Let $L$ and $\LL$ be as in Lemma 4.3. Consider the Hamiltonian cobordisms
$\beta_H$ associated to (positively) adimissible Hamiltonian $H$. Then 
two such Hamiltonian cobordisms are Hamiltonian isotopic to each other in
$P \times T^*[0,1]$ by an isotopy that is compactly supported in
$P \times T^*(0,1)$.
\end{lem}

\medskip

We now study the condition (2) from Proposition 4.2
 that {\it the Hamiltonian cobordism $\beta_H$
is transverse to the level sets of $\widetilde \varphi = \varphi + {1\over 2}
(s^2 +a^2)$}. A typical example of such Hamiltonians arise in the
following way: Let $\partial P = M$ with its induced contact structure
and $L_0 \subset P$ be a proper Lagrangian submanifold
with its boundary $R_0 \subset M$. $R_0$ is a compact Legendrian
submanifold of $M_0$. Consider a Hamiltonian isotopy of $L_0$
which extends a Legendrian isotopy of $R_0 \subset M$. We choose
Hamiltonians which restrict to {\it contact Hamiltonians} (see [A2]
for the definition) on the
collar $(1-\epsilon, 1] \times M$ of $\partial P$, i.e., satisfies
\begin{equation}
H(cm, t) = cH(m,t) \quad \mbox{\rm for } \, m \in M, \, c\in \R^+
\end{equation}
on the symplectic cone attached to $\partial P$.
\medskip

\begin{lem}
Let $P$ be Weinstein at infinity and $\varphi$ be an
exhausting pluri-subharmonic function which is super-quadratic over the radial
coordinate. Suppose that $H$ satisfy (15) and
that $L_0$ is transverse to the level sets of $\varphi$ at infinity. Then
the induced Hamiltonian cobordism $\beta_H \subset P \times T^*\R$ of $L_0$
is transverse to the level sets of $\widetilde \varphi$ at infinity.
\end{lem}

\begin{proof}
Since we extend $\beta_H$ so that $a = 0$ outside $0\leq s \leq 1$
which is obviously transverse, it is enough to check the transversality over
$0\leq s \leq 1$. In this region, we may consider the function $\varphi + {1\over
2} a^2$ in place of $\widetilde \varphi$. Recalling
$$
\beta_H=\{ (x,s,a)\in P\times  T^*\R\mid x\in L_s,\; a=-H(x,s),\; 0\le s\le
1\},
$$
it is easy to check that the tangent space of $\beta_H$
at $(x,s,-H(x,s))$ is spanned by the vectors
$$
\vec{v} - H(x,s)dH(\vec{v}){\partial \over
\partial a} + c\Big({\partial \over \partial s}- {\partial H \over \partial
s}{\partial \over \partial a}\Big)
$$
where $d$ is the differential for $x$, $\vec{v}
\in T_xL_s$ and $c \in \R$. Applying this vector to $\varphi + {1\over 2}a^2$, the
non-transversal points are characterized by the equation
\begin{eqnarray}
 &d\varphi(x)\big|_{L_s} + H(x,s)  dH(x,s)\big|_{L_s}=0 \nonumber \\
 &H(x,s)  {\partial H \over \partial s} = 0, \quad a = -H(x,s)
\end{eqnarray}
on the collar or on the symplectic cone attached to $\partial P$. Since
$d\varphi \neq 0$, $H(x,s) \neq 0$ on the cone. On the other hand, since
the growth of $H$ is linear and the growth of $\varphi$
is super-quadratic over the radial coordinate, the first equation of (16)
cannot hold at infinity
in $P \times T^*[0,1]$. This finishes the proof. 
\end{proof}

\medskip

We now apply Theorem 2.1 to the case $L_0 = \alpha_L$ and $L_1 = \beta_H$
to obtain the $C^0$-estimate for (10). Once the crucial $C^0$-estimate
is obtained, the standard arguments in the Floer theory
prove the following proposition.
\medskip

\begin{thm} 
Let $\partial_0:CF(L,L_0)\rightarrow CF(L,L_0)$
and $\partial_1 :CF(L,L_1)\rightarrow CF(L,L_1)$ be the Floer boundary maps on
$(P,\omega )$. Suppose $L,  {\mathcal L}$ satisfy the properties
required in Definition 4.2 and let $H$ a Hamiltoian generating $\mathcal L$
and satisfying (14) and (15).
 Let $h_{\beta_H ,\e}:CF(L,L_0:\alpha_\e )\rightarrow
CF(L,L_1:\alpha_\e )$ be the map defined in (11).  Then the identity (12)
holds and so $h_{\beta_H ,\e}$'s induce a homomorphism, as $\e \to 0$,
$$
h_{\beta_H ,\e} :HF(L,L_0:J)\rightarrow HF(L,L_1:J).
$$
Furthermore, this
homomorphism is independent of the approximations $\alpha_\e$ and of the choice of
$H$. We denote the common homomorphim by
\begin{equation}
h_{\mathcal L}: HF(L,L_0) \to HF(L, L_1)
\end{equation}
\end{thm}
\begin{proof} 
Under the hypotheses given in the statement, it follows that
$$
\widetilde y\in\partial\beta\cap\alpha_{L,\e}=L\cap L_0\times\{
(0,0)\}\amalg L\cap L_1\times \{(1,0)\}.
$$
Once we have this, the standard argument in the Floer theory proves
the chain property (4.7).

To prove the independence of $h_{\beta_H ,\e}$ on $H$ and $\e >0$, it will be enough
to prove that the family of approximations $\{\alpha_\e\}_{\e >0}$ and the
change of $H$'s satisfying the condition above can be realized
by compactly supported Hamiltonian deformations of one another among them.  But this
follows from the construction of  the approximation $\Upsilon_\e$ of
$Ch(1_{[0,1)})$. We refer to [KO1,2] for the details of this limiting argument.
\end{proof}

In the next sections, we will prove that our chain map associated
to a Hamiltonian isotopy $\mathcal L$ is natural and becomes an isomorphism.

\medskip

\begin{rem} 
(1)~ This will fill the gap present in the construction of
the chain isomorphism used in [O2], which the author overlooked in
applying the strong maximum principle to get the $C^0$-estimate
for the continuity equation (2) or (3).
This $C^0$-estimate and the isomorphism were crucial in the proof of
continuity of the invariants $S \mapsto \rho(H,S)$ under the isotopy of
submanifolds $S$ (see the proof of Proposition 6.5 [O2]).
\smallskip

(2)~ The Hamiltonian isotopies considered in Lemma 4.5 includes
the {\it positive Lagrangian isotopy} of {\it $\theta$-exact Lagrangian
subamnifolds} considered in [KhSe]. In particular, we have provided
the recipe of curing the ``weakness'' mentioned therein
in that our construction provides a canonical isomorphism to
Lemma 5.11 [KhSe] that was missing therein.
\end{rem}

\section{Composition rule}

In this section, we will prove the following composition rule,

\begin{equation} 
h_{\beta_0\# \beta_1} =h_{\beta_1}\circ h_{\beta_0}
\end{equation} 
where $L_0
{\displaystyle{\mathrel{\mathop{\sim}_{\beta_0}}}} L_1$,
$L_1{\displaystyle{\mathrel{\mathop{\sim}_{\beta_1}}}} L_2$ and $\beta_0\#\beta_1$
denotes the obvious composition of Lagrangian cobordisms $\beta_0$ and $\beta_1$.

We examine how the Lagrangian boundary conditions are involved.  At $t=1$, we can
just take a small  perturbation of the elongated $\beta_0\#\beta_1$.  At $t=0$, we
need to describe some approximation result for
$$
L\times Ch(1_{[0,1)})\cup L\times Ch(1_{[1,2)})=L\times
(Ch(1_{[0,1)})\cup Ch(1_{[1,2)})).
$$

\begin{figure}[htb]
\includegraphics{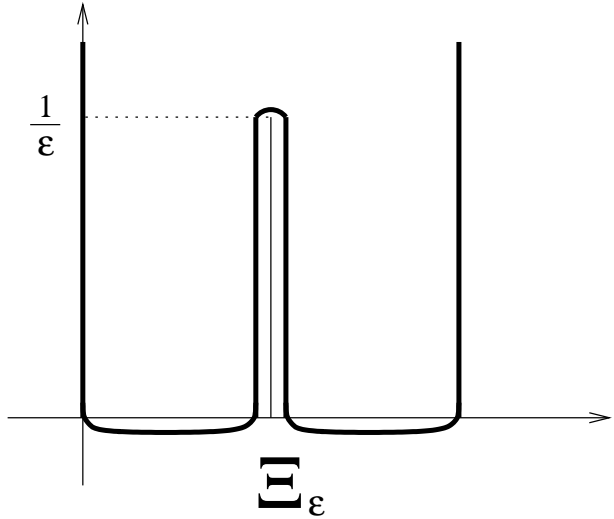}
\caption{Approximation of the cycle $Ch(1_{[0,1)}) \cup Ch(1_{[1,2)})$}
\end{figure}

\n by a family of Lagrangian submanifolds $\{\Xi_\e\}_{0<\e <1}$ as drawn
above (See [KO2] for many illustrations of this approximation argument).
First, we remark that for any given compact subset of $\e$ in
$(0,1)$, the corresponding Lagrangian submanifolds $L\times\Xi_\e$ are
deformations to one another via {\it compactly supported} Hamiltonian isotopies
$T^*(M\times \R)$.
Then some modification of standard gluing arguments can be applied to prove
the following analytical result (See [KO2] for some relevant discussion).

\medskip

\begin{thm} 
There exists sufficiently small $\e >0$ such that we
have gluing diffeomeorphisms
$$ \MM (\Jbar, \alpha_\e ,\beta_0)\times\MM
(\Jbar,\alpha_\e ,\beta_1)\rightarrow \MM (\Jbar, L\times\Xi_\e,
\beta_0\#\beta_1)
$$
after a modification of the cobordism $\beta_0 \# \beta_1$ near $s =1$
as described in the proof of Lemma 4.5.
In particular, we have the identity,
\begin{equation} 
h_{\beta_1,\e}\circ
h_{\beta_0,\e}=h_{L\times\Xi_\e}:HF_*(L,L_0)\rightarrow HF_*(L,L_2).
\end{equation} 
\end{thm}

\begin{proof} 
We will be sketchy in the proof because similar gluing arguments
have been used many times in the literature by now.

Note that $\beta_0 \# \beta_1$ is again a Hamiltonian cobordism. We choose
a Hamiltonian that is postively admissible 
to $\beta_0 \# \beta_1$. In fact,
by adding a bump function supported in a neighborhood of the hypersurface $s =1$,
we can make the corresponding Hamiltonian $H$ so that 
$\hbox{\rm Graph }H$ is ``above'' $\Xi_\e$ as in Figure 4. We glue each given pair
$u_0 \in \MM(\widetilde J,\alpha_\e,\beta_0)$  and
$u_1 \in \MM(\widetilde J,\alpha_\e,\beta_1)$ with
the obvious holomorphic strip in the middle. This gluing
is possible, as long as $\e$ is sufficiently small and so ${1\over \e}$
is sufficiently large and the strip is sufficiently narrow. This finishes the proof.
\end{proof}

\begin{figure}[htb]
\includegraphics{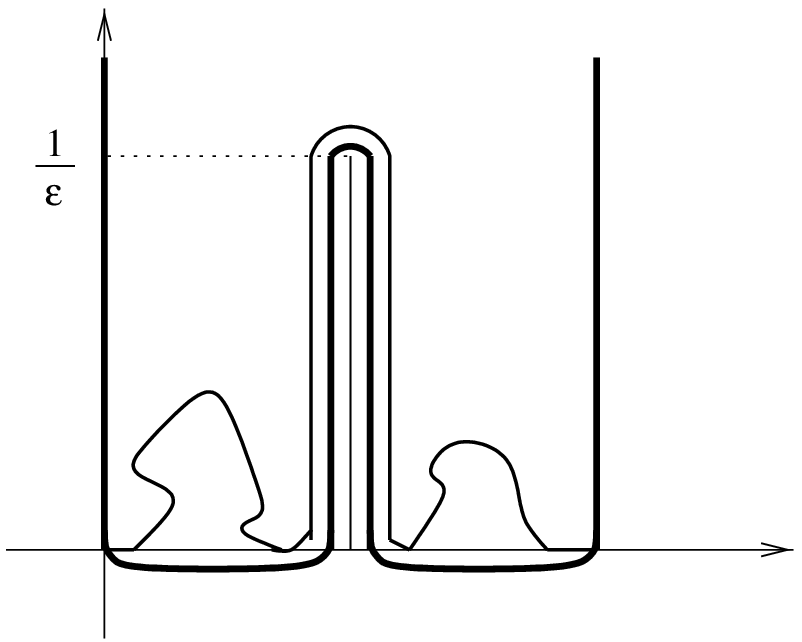}
\caption{}
\end{figure}

After this crucial analytical step, we use the fact, which can be easily checked,
 that the family $\{ L\times
\Xi_\e\}_{0<\e <1}$ are Hamiltonian deformations to one another via {\it compactly
supported} Hamiltonian isotopy.  Therefore we can apply the standard continuation
argument in the Floer theory to show that the homomorphisms
$$ 
h_{L\times\Xi_\e}
:HF_*(L,L_0)\rightarrow HF_*(L,L_1)
$$
are independent of $\e >0$.  Since we have
\begin{equation}
h_{(\beta_0\#\beta_1)}=\lim_{\e\to 0} h_{(L\times \Xi_\e)},
\end{equation}
we have finished proof of (18) combining (19) and (20).
\medskip

\section{Trivial cobordism}

In this section we prove the following theorem. This is the place where the
power of choosing $\alpha_L$ as we do for the boundary condition at $t=0$
becomes manifest.

\medskip

\begin{thm} 
Consider the trivial product cobordism
$$
\beta_0 = L_0 \times [0,1] \times \{0\} \subset P \times T^*\R.
$$
Then the induced homomorphism
$h_{\beta_0}: HF_*(L, L_0) \to HF_*(L, L_0)$ is the identity homomorphim.
\end{thm}

\begin{proof} 
Recall $\alpha_\e = L \times \Xi$. We study the equation
\begin{equation}
\left\{
\begin{array}{l}
{\partial\widetilde u\over\partial\tau} +\widetilde
J\Big( {\partial\widetilde u\over\partial t}\Big) =0\cr 
\widetilde u(\tau ,0)\in \alpha_\e,\; \widetilde u(\tau ,1)\in\beta_0
\end{array}\right. 
\end{equation}
Since
$\widetilde J = J\oplus i$, $\alpha_\e = L \times \Upsilon_\e$ and $\beta_0 = L_0 \times
o_{\R}$ all split, (21) splits into
\begin{equation}
\left\{ 
\begin{array}{l}{\partial u\over\partial\tau}+J{\partial u\over\partial t}=0\cr 
u(\tau ,0)\in L,\;\; u(\tau ,1)\in L_0 
\end{array}\right.
\end{equation}
and
\begin{equation}
\left\{ 
\begin{array}{l}
{\overline \partial v} = 0  \cr 
v(\tau,0) \in \Upsilon_\e, \;\; v(\tau,1) \in o_\R. \end{array}
\right.
\end{equation}
Noting that (23) has the unique solution with index 1 (up to translations) and
with the asymptotic condition
\begin{equation}
v(-\infty) = (0,0), \;\; v(\infty) = (1,0),
\end{equation}
solutions $\widetilde u \in
\MM_\e(\xbar,\ybar: \Jbar, \alpha_\e, \beta_0) $ of (21) with index 1 consist of
the pairs $(u,v)$ such that $u$ is a solution  $u$ of (22) with index 0 and $v$
is the unique solution of (23) satisfying (24). In particular, $u$ must be
constant. Hence we have proven that
\begin{eqnarray*}
n_\e(\xbar,\xbar:\beta_0) & = & 1  \quad \mbox{\rm for all }\, x \in L\cap
L_0 \\ 
n_\e(\xbar,\ybar:\beta_0) & = & 0  \quad  \mbox{\rm if } \, y\neq x 
\end{eqnarray*}
which in turn implies that the chain map
$$
h_{\beta_0,\e}: CF(L,L_0:J, \alpha_\e) \to CF(L,L_0:J, \alpha_\e)
$$
becomes the identity map. This finishes the proof. 
\end{proof}

One immediate corollary of (18) and Theorem 6.1 is the following

\medskip

\begin{thm} 
Let $H$ be a positively admissible Hamiltonian to $(L, \LL)$,
and $\beta_H$ be the Hamiltonian cobordism obtained from
the Hamiltonian isotopy $\phi^s_H(L_0)$ from $L_0$ to $L_1=\phi^1_H(L_0)$.  Then
the homomorphism $$ h_{\beta_H} :HF_*(L,L_0)\rightarrow HF_*(L,L_1) $$ is an
isomorphism. Hence $h_{ \mathcal L}: HF_*(L,L_0) \to HF_*(L,L_1)$ is an
isomorphism. 
\end{thm}

\begin{proof}
We compose $\beta_H$ with $\beta_{\overline H}$ where

$${\overline H}(x,s) :=-H(\phi^s_H(x), s)$$ which generates the isotopy $\{
(\phi^s_H)^{-1}(L_1)\}$.  It is immediate to check that the composition
$\beta_H\#\beta_{\overline H}$ is Hamiltonian isotopic to the product cobordism
between $L_0$ and $L_0$ via {\it compactly supported} Hamiltonian isotopy
$P \times T^*\R$.
Therefore we can apply the standard procedure of using (4) to prove the
construction of chain isomorphisms between the cases of the identity
cobordism and $\beta_H\#\beta_{\overline H}$. This proves the theorem.
\end{proof}

\section{Intersection of conormal bundles}

In this section, we apply our extended Floer theory to the special case of
conormal bundles $\nu^*S_1, \nu^*S_2$ of two smooth submanifolds
$S_1, S_2 \subset M$.  We would like to compute $HF_*(\nu^*S_1, \nu^*S_2)$,
when $S_1$ is transverse to $S_2$.

First we note that the intersection of conormals
$$
\nu^*S_1 \cap \nu^*S_2 = o_{S_1\cap S_2}
$$
is compact and the following types of deformations or compositions of them
leave the intersection set compact:

\smallskip

(1) $\phi_t$ are compactly supported, or

(2) they are homogeneous symplectomorphisms (at infinity), i.e., it is generated
by the Hamiltonian of the form $(q,p) \mapsto \langle p, X_t(q) \rangle$
such that {\it $S_1$ is transverse to $f_t(S_2)$ for all $t$} where $f_t:M\to M$
is the flow of $X_t$, or

(3) they are the fiberwise translations by $t df$ where $f$ is a smooth
function defined on the base $M$.

\smallskip
One can easily check that any two such $\Phi= \{\phi_t\}$ can be connected by
one parameter family $\{\Phi^s\}_{0\leq s\leq 1}$ such that intersections of $\nu^*S_1$ and
$\phi_t^s(S_2)$ remain to be compact. Therefore it follows from  the discussions
in the previous sections that there exist a canonical chain isomorphism
$$
h:CF(\nu^*S_1,\phi_1(\nu^*S_2)) \to CF(\nu^*S_1, \phi_2(\nu^*S_2))
$$
where, for example, $\Phi_i=\{\phi^t_i\}_{0\leq t\leq 1}$ for $i = 1,\, 2$
is a Hamiltonian isotopy of $T^*M$ of the above types or a composition of them.
Therefore this induces the canonical isomorphism
$$
h:HF(\nu^*S_1,\phi_1(\nu^*S_2)) \to HF(\nu^*S_1, \phi_2(\nu^*S_2)).
$$
We denote by $HF(\nu^*S_1, \nu^*S_2)$ the common group.

The existence of such isomrophisms for the first two cases is immediate
from the discussions in the previous sections.  The case (3) follows
since we can easily check that the corresponding Hamiltonian cobordism
satisfies the hypotheses (1) and (2) from Proposition 4.2.

To compute $HF_*(\nu^*S_1, \nu^*S_2)$, we deform $\nu^*S_2$ to
$\phi_f(\nu^*S_2)$ where $\phi_f$ is
the fiberwise translations by $df$, where the function $f$ on M
will be suitably chosen.  Under preliminary perturbation of the class (2)
of $S_2$ in $M$, we may assume that $S_1$ is transverse to  $S_2$ 
and so $S_1\cap S_2$
is a smooth submanifold. We choose a smooth Morse function
$\widetilde f: S_1\cap S_2 \to \R$ and extend it to $M$, first
quadratically to a tubular neighborhood and then suitably cutting
off outside the neighborhood (See [Pz] or [O2]). We denote the extension
by $f:M \to \R$.
\medskip

\begin{prop}
Let $f$ and $\phi_f$ described as above. Then
we have
\smallskip

$(1)$ $\nu^*S_1$ is transverse to $\phi_f(\nu^*S_2)$

$(2)$ $\nu^*S_1 \cap \phi_f(\nu^*S_2)$ is finite and all lie in the zero
section of $T^*M$.
\end{prop}

\begin{proof} 
We first prove (2). Let $\alpha_1 \in \nu_q^*S_1$. If
$\alpha_1 \in \nu_q^*S_1\cap \phi_f(\nu^*S_2)$, then we should have
\begin{equation}
\alpha_1 =  \alpha_2 + df(q)
\end{equation}
for some $\alpha_2 \in \nu^*_qS_2$. Since $\alpha_1 \in \nu^*_qS_1$,
we must have
\begin{equation}
\alpha_2|_{TS_1} = - df(q)|_{TS_1}
\end{equation}
Since $S_1$ is transverse to $S_2$ and $\alpha_2 \in \nu^*_qS_2$, (7.2)
uniquely determines $\alpha_2$ and so $\alpha_1$. It remains to
show that $\alpha_1 = 0$. To show this, it is enough to prove that
$$
\alpha_1|_{TS_2/T(S_1 \cap S_2)}\equiv 0
$$
because $\alpha_1|_{TS_1}\equiv 0$. On the other hand, this follows from
(25) noting that $\alpha_2|_{TS_2} \equiv 0$ and that we have extended
the Morse function $\widetilde f$ on $S_1\cap S_2$ quadratically to
its tubular neighborhood and so $df(q)|_{TM/T(S_1\cap S_2)} \equiv 0$.
This finishes the proof of (2). Once we prove this, (1) immediately
follows from transversality of intersections of $S_1$ and $S_2$.
\end{proof}
\medskip

With Propostions 7.1 in our hand, we can repeat the computations from [F2], [Pz]
or [O2] to construct one to one correpondence between the moduli space
$\MM(J, \nu^*S_1, \phi_f(\nu^*S_2))$ of Floer's trajectories and the moduli
space $\MM^{Morse}(f; S_1 \cap S_2)$ for a suitably chosen almost complex
structure $J$ (see [F2], [Pz] for the relevant arguments in a different context).
Combining these and construction of orientation of the Floer moduli space
from [Oh2], we have proved the following

\medskip

\begin{thm}
Let $S_1, S_2 \subset M$ be a compact smooth submanifolds
and $\nu^*S_1, \, \nu^*S_2$ be their conormal bundles. Then there exists a
canonical chain isomorphim
$$
C^{Morse}(f;S_1\cap S_2) \to CF(\nu^*S_1, \phi_f(\nu^*S_2))
$$
which induces an isomorphism
$$
h: H_*(S_1\cap S_2;\Z_2) \to HF(\nu^*S_1, \phi_f(\nu^*S_2))
\simeq HF(\nu^*S_1,\nu^*S_2)
$$
in $\Z_2$-coefficients in general.
When $S_1, S_2$ and $M$ are oriented, then this isomorphism holds in
$\Z$-coefficients. 
\end{thm}

\medskip

This combined with the invariance property of the Floer homology
under the Hamiltonian isotopy of the types, e.g., (1), (2) and (3) above,
immediately gives rise to the following intersection theorem.

\medskip

\begin{cor} Let $S_1, \, S_2$ be as before.
Suppose $\phi$ is a Hamiltonian diffeomorphism
on $T^*M$ of the types above or a composition of them. Then
$$
\# (\nu^*S_1 \cap \phi(\nu^*S_2)) \geq \hbox{\rm rank } H_*(S_1\cap S_2)
$$
provided $\nu^*S_1$ is transverse to $\phi(\nu^*S_2)$.
Here $H_*(S_1\cap S_2)$ is in $\Z$-coefficients in the
oriented case and in $\Z_2$-coefficients in general.
\end{cor}
\medskip

We would like to compare Theorem 7.2 with the conjecture stated
in the end of [GM]. It would be very interesting to generalize the
construction in [KO1,2] to the general stratified case to give
a precise meaning of the statement of the conjecture [GM].

\section{Further discussions}

In [Po], Polterovich introduced the notion of {\it Lagrangian pseudo-isotopy}
and in [C], Chekanov introduced that of (connected) {\it monotone Lagrangian
cobordism}. If we restrict to the case of
monotone Lagrangian submanifolds for which the Floer homology can be
easily constructed without any sophisticated machinery,
the construction we have carried out in the previous sections
also applies to the monotone Lagrangian cobordism, in particular to
the Lagrangian pseudo-isotopy. Therefore we have proved that for any monotone
Lagrangian cobordism $\beta$ from $L_0$ and $L_1$, there exists a natural
homomorphism
$$
h_{\beta}: HF(L, L_0) \to HF(L, L_1).
$$
For more complicated cobordism, we do not expect such homomorphims but
expect only some ``correspondences''.

In fact, this construction works for the case of Lagrangian pseudo-isotopy as long
as the Floer homology $HF(L,L_0)$ for the given Lagrangian submanifold
$L$ and $L_0$ can be constructed (We refer to [FOOO] for the most general
construction of Floer homology upto now). Unlike the case of Hamiltonian isotopy,
the corresponding chain map is not expected to be an isomorphism and so can
provide an obstruction to Lagrangian pseudo-isotopy being a Hamiltonian isotopy.
It would be an extremely interesting problem to find a nontrivial
Lagrangian pseuo-isotopy, when there is.

One very interesting problem is to study the change of $HF(L,L^\prime)$ when
the isotopy $\{L_t\}_{0\leq t \leq 1}$ of $L^\prime$ undergoes the process of losing
the intersections to infinity. A model case to study will be the one
of symplectic manifolds with contact type boundary and its proper
Lagrangian submanifolds. In this case, the corresponding family of
boundary Legendrian submanifolds will have intersections at a finite
number of $t$'s in (0,1) with $L$. In particular, it would be interesting
to describe the change of $HF(\nu^*S_1, \nu^*S_2^t)$ at the time $t_0$, where
the intersection pattern of $S_1\cap S_2^t$ changes  through a degenerate intersection.
This will be a subject of future study.

\bigskip

{\bf Acknowledgements:}
The idea of the present paper was first presented in the Symplectic Geometry
Workshop in Warwick University in the summer of 1998. We thank D. Salamon for
the invitation and for some useful discussions.
We also thank K. Fukaya for some helpful discussions
during our visit of RIMS in the fall of 1999,
L. Polterovich for drawing our attention to the paper [C] and
K. Hori for a very inspiring lecture in KIAS on the result from [HIV].


\begin{thebibliography}{FOOO}
\bibitem[A1]{} Arnold, V.\ I., {\it Lagrange and Legendre cobordism, I},
Funkt.\ Anal.\ Ego Prilozh.  {\bf 14} (1980), 1--13. (Funct. Anal. Appl. (1981),
167--177).

\bibitem[A2]{} Arnold, V.\ I., Mathematical Methods of Classical Mechanics,
Springer-Verlag, New York, 1978.

\bibitem[C]{} Chekanov, Yu. V., {\it Lagrangian embeddings and Lagrangian cobordisms},
Amer. Math. Soc. Transl. (2) vol 180, 1997, pp 13-23, edited by Kovansky,
Varchenko \& Vassiliev.

\bibitem[EG1]{} Eliashberg, Y. and Gromov, M., {\it Convex symplectic manifolds,}
in Several Complex Variables and Complex Geometry (I. Bedford, etc eds.), Proc.
Sympos. Pure Math. 52, Part 2, AMS, Providence, RI, 1991, pp 135--162.

\bibitem[EG2]{} Eliashberg, Y. and Gromov, M., {\it Lagrangian intersections theory:
finite-dimensional appropach,} in Geometry of Differential Equations, 1998, pp 27-118.

\bibitem[EHS]{} Eliashberg, Y., Hofer, H. and Salamon, D., {\it Lagrangian
intersections in contact geometry,}  Geom. Funct. Anal. {\bf 5} (1995),
244--269. 

\bibitem[F1]{} Floer, A., {\it Morse theory for Lagrangian intersections,} J.
Differ. Geom. {\bf 28} (1988), 513--547.

\bibitem[F2]{} Floer, A., {\it Witten's complex and infinite dimensional Mores
theory,} J. Differ. Geom. {\bf 30} (1989), 207--221.

\bibitem[F3]{} Floer, A., {\it Symplectic fixed points and holomorphic spheres},
Commun.\ Math.\ Phys.\ {\bf 120} (1989), 575--611.

\bibitem[FOOO]i{} Fukaya, K., Y.-G. Oh, H. Ohta and K. Ono,
{\it Lagrangian intersection Floer theory - anomaly and obstruction -,}
(to appear).

\bibitem[GM]{} Grinberg, M. and MacPherson, R., {\it Euler characteristics and
Lagrangian intersections,} in IAS/Park City Math. Series, vol 7, 1999, edited by
Eliashberg, Y. and Traynor, L. AMS, 1999, pp 267-291.

\bibitem[HIV]{} Hori, K., Iqbal, A. and Vafa, C., {\it $D$-branes and mirror symmetry,}
preprint 2000/ hep-th/0005247.

\bibitem[KO1]{} Kasturirangan, R., and Oh, Y.-G, {\it Floer homology of open
subsets and a relative version of Arnold's conjecture,} Math. Z., (to appear).

\bibitem[KO2]{} Kasturirangan, R., and Oh, Y.-G., {\it Quantization of
Eilenberg-Steenrod axioms via Fary functors,} submitted.

\bibitem[KaSc]{} Kashiwara, M. and Schapira, P., Sheaves on Manifolds, A Series of
Comp. Studies in Math. vol 292, Springer-Verlag, New York, 1990.

\bibitem[KhSe]{} Khovanov, M. and Seidel, P., {\it Quivers, Floer cohomology, and
braid group actions,} preprint 2000/ math.QA/ 0006056.

\bibitem[O1]{} Oh, Y.-G., {\it Floer cohomology of Lagrangian intersections and
pseudo-holomorphic disks}, I, Comm.\ Pure Appl.\ Math.\ {\bf 46} (1993),
949--994; Addenda, ibid.  {\bf 48} (1995), 1299--1302.
\smallskip

\bibitem[O2]{} Oh, Y.-G., {\it Symplectic topology as the geometry of action
funtional I}, J. Differ. Geom. {\bf 46} (1997), 499--577.

\bibitem[Po]{} Polterovich, L., {\it Symplectic displacement energy for
Lagrangian submanifolds,} Ergodic Theory and Dynam. Sys., 13 (1993),
357--367.

\bibitem[Pz]{} Po\'zniak, {\it Floer homology, Novikov rings and clean intersections,}
Thesis in University of Warwick, 1994.

\end{thebibliography}
\end{document}